\documentclass[11pt]{amsart}
\usepackage{amssymb}

\DeclareMathOperator{\dom}{Dom} 

\newcommand{\dbar}{\ensuremath{\overline\partial}}
\newcommand{\dbarstar}{\ensuremath{\overline\partial^*}}
\newcommand{\C}{\ensuremath{\mathbb{C}}}
\newcommand{\R}{\ensuremath{\mathbb{R}}}

\makeatletter
\newcommand{\sumprime}{\if@display\sideset{}{'}\sum%
            \else\sum'\fi}
\makeatother

\begin{document}

\numberwithin{equation}{section}

% define theorem environments
\newtheorem{theorem}{Theorem}[section]
\newtheorem{proposition}[theorem]{Proposition}
\newtheorem{conjecture}[theorem]{Conjecture}
\def\theconjecture{\unskip}
\newtheorem{corollary}[theorem]{Corollary}
\newtheorem{lemma}[theorem]{Lemma}
\newtheorem{observation}[theorem]{Observation}
\theoremstyle{definition}
\newtheorem{definition}{Definition}
\numberwithin{definition}{section}
\newtheorem{remark}{Remark}
\def\theremark{\unskip}
\newtheorem{question}{Question}
\def\thequestion{\unskip}
\newtheorem{example}{Example}
\def\theexample{\unskip}
\newtheorem{problem}{Problem}

\def\vvv{\ensuremath{\mid\!\mid\!\mid}}
\def\intprod{\mathbin{\lr54}}
\def\reals{{\mathbb R}}
\def\integers{{\mathbb Z}}
\def\N{{\mathbb N}}
\def\complex{{\mathbb C}\/}
\def\distance{\operatorname{distance}\,}
\def\spec{\operatorname{spec}\,}
\def\interior{\operatorname{int}\,}
\def\trace{\operatorname{tr}\,}
\def\cl{\operatorname{cl}\,}
\def\essspec{\operatorname{esspec}\,}
\def\range{\operatorname{\mathcal R}\,}
\def\kernel{\operatorname{\mathcal N}\,}
\def\linearspan{\operatorname{span}\,}
\def\lip{\operatorname{Lip}\,}
\def\sgn{\operatorname{sgn}\,}
\def\Z{ {\mathbb Z} }
\def\e{\varepsilon}
\def\p{\partial}
\def\rp{{ ^{-1} }}
\def\Re{\operatorname{Re\,} }
\def\Im{\operatorname{Im\,} }
\def\dbarb{\bar\partial_b}
\def\eps{\varepsilon}

\def\Hs{{\mathcal H}}
\def\E{{\mathcal E}}
\def\scriptu{{\mathcal U}}
\def\scriptr{{\mathcal R}}
\def\scripta{{\mathcal A}}
\def\scriptc{{\mathcal C}}
\def\scriptd{{\mathcal D}}
\def\scripti{{\mathcal I}}
\def\scriptk{{\mathcal K}}
\def\scripth{{\mathcal H}}
\def\scriptm{{\mathcal M}}
\def\scriptn{{\mathcal N}}
\def\scripte{{\mathcal E}}
\def\scriptt{{\mathcal T}}
\def\scriptr{{\mathcal R}}
\def\scripts{{\mathcal S}}
\def\scriptb{{\mathcal B}}
\def\scriptf{{\mathcal F}}
\def\scriptg{{\mathcal G}}
\def\scriptl{{\mathcal L}}
\def\scripto{{\mathfrak o}}
\def\scriptv{{\mathcal V}}
\def\frakg{{\mathfrak g}}
\def\frakG{{\mathfrak G}}

\def\ov{\overline}
%\date {November, 2004.  %Print \today}

\author{Siqi Fu}
\thanks
{This research was supported in part by  an NSF grant.}
\address{Department of Mathematical Sciences,
Rutgers University-Camden, Camden, NJ 08102}
\email{sfu@camden.rutgers.edu}
\title[]  % runningheader
{Spectrum of the $\dbar$-Neumann Laplacian on polydiscs} \maketitle

\begin{abstract} The spectrum of the $\dbar$-Neumann Laplacian on a polydisc in
$\C^n$ is explicitly computed.  The calculation exhibits that the
spectrum consists of eigenvalues, some of which, in particular the
smallest ones, are of infinite multiplicity.
\end{abstract}

\section{Introduction}\label{intro}

The $\dbar$-Neumann Laplacian $\square_q$ on a bounded domain
$\Omega$ in $\C^n$ is (a constant multiple of) the usual Laplacian
acting diagonally on $(0, q)$-forms subjected to the non-coercive
$\dbar$-Neumann boundary conditions. It is a densely defined,
non-negative, and self-adjoint operator. As such, its spectrum is
a non-empty closed subset of the non-negative real axis.  Unlike
the usual Dirichlet Laplacian, its spectrum needs not be purely
discrete. (See \cite{FuStraube01} for a discussion on related
subjects.) Spectral behavior of the $\dbar$-Neumann Laplacian is
more sensitive to the boundary geometry of the domain than the
Dirichlet/Neumann Laplacians. (See \cite{Fu05a, Fu05b} and
references therein for related discussions.)

Spectral behavior of the $\dbar$-Neumann Laplacian on special
domains often serves as a model for the general theory.  One
certainly cannot expect to explicitly calculate the spectrum for
wide classes of domains.  The spectrum for the ball and annulus
was explicitly computed by Folland \cite{Folland72}. In this note,
we compute the spectrum for the polydiscs.  Our computation
exhibits that the spectrum of the $\dbar$-Neumann Laplacian on a
polydisc consists of eigenvalues, some of which, in particular the
smallest ones, are of infinite multiplicity. That the essential
spectrum of the $\dbar$-Neumann Laplacian is non-empty is
consistent with, in fact, equivalent to, the well-known fact that
the $\dbar$-Neumann operator (the inverse of the $\dbar$-Neumann
Laplacian) is non-compact (e.g., \cite{Krantz88}). It is
noteworthy that for a polydisc, the bottom of the spectrum is
always in the essential spectrum---a phenomenon not stipulated in
the general operator theory.

\section{Preliminaries}\label{prelim}

We first recall the setup for the $\dbar$-Neumann Laplacian. We
refer the reader to \cite{FollandKohn72, ChenShaw99} for an in
depth treatise of the $\dbar$-Neumann problem.

Let $\Omega$ be a bounded domain in $\C^n$. For $1\le q\le n$, let
$L^2_{(0, q)}(\Omega)$ denote the space of $(0, q)$-forms with
square integrable coefficients and with the standard Euclidean
inner product whose norm is given by
\[
\|\sumprime a_{J}d\bar z_J\|^2= \sumprime \int_\Omega |a_{J}|^2
dV(z),
\]
where the prime indicates the summation over strictly increasing
$q$-tuples $J$. (We consider $a_J$ to be defined on all
$q$-tuples, antisymmetric with respect to $J$.) For $0\le q\le
n-1$, let $\dbar_q\colon L^2_{(0, q)}(\Omega)\to L^2_{(0,
q+1)}(\Omega)$ be the $\dbar$-operator defined in the sense of
distribution. This is a closed and densely defined operator. Let
$\dbarstar_q$ be its adjoint. Then $\dbarstar_q$ is also a closed
and densely defined operator with domain
\[
\dom(\dbarstar_q)=\{u\in L^2_{(0, q+1)}(\Omega) \mid \exists C>0
\text{ such that } |\langle u, \dbar v\rangle|\le C \|v\|, \forall
v\in\dom(\dbar_q)\}.
\]
For $1\le q\le n-1$, let
\[
Q_q(u, v)=(\dbar_q u, \dbar_q v)+(\dbarstar_{q-1} u,
\dbarstar_{q-1} v)
\]
be the sesquilinear form on $L^2_{(0, q)}(\Omega)$ with
$\dom(Q_q)=\dom(\dbar_q)\cap\dom(\dbarstar_{q-1})$. It is evident
that $Q_q$ is non-negative, densely defined, and closed. The
$\dbar$-Neumann Laplacian
$\square_q=\dbarstar_q\dbar_q+\dbar_{q-1}\dbarstar_{q-1}\colon
L^2_{(0, q)}(\Omega)\to L^2_{(0, q)}(\Omega)$ is the associated
self-adjoint operator with domain $\dom(\square_q)=$
\[
\{u\in L^2_{(0, q)}(\Omega) \mid u\in
\dom(\dbar_q)\cap\dom(\dbarstar_{q-1}), \dbar u\in
\dom(\dbarstar_q), \dbarstar_{q-1} u\in \dom(\dbar_{q-1})\}.
\]

For the reader's convenience, we also briefly review relevant
facts of the Bessel functions. Extensive treatment of the Bessel
functions can be found, for example, in \cite{Watson48}.  The
Bessel functions of integer orders are given via the following
Laurent expansion:
\begin{equation}\label{laurent}
e^{\frac{z}{2}(t-\frac{1}{t})}=\sum_{m=-\infty}^\infty t^m J_m(z).
\end{equation}
Evidently, $J_{-m}(z)=(-1)^m J_m(z)$ and when $m\ge 0$,
\[
J_m(z)=\sum_{l=0}^\infty \frac{(-1)^l(z/2)^{2l+m}}{l!(l+m)!}.
\]
Hence $J_m(z)$ is an entire function with zero of order $|m|$ at
the origin. By differentiating both sides of \eqref{laurent} with
respect to $t$ and with respect to $z$, we have the recurrence
formulas:
\begin{equation}\label{rec-1}
mJ_m(z)=\frac{z}{2}(J_{m+1}(z)+J_{m-1}(z)), \qquad
J'_m(z)=\frac{1}{2}(J_{m-1}(z)-J_{m+1}(z)).
\end{equation}
Therefore,
\begin{equation}\label{rec-2}
zJ_{m-1}(z)=zJ'_m(z)+mJ_m(z), \qquad
zJ_{m+1}(z)=-zJ'_n(z)+mJ_m(z).
\end{equation}
It follows that $J_m(z)$ satisfies the Bessel equation:
\begin{equation}\label{bessel-eq}
J''_m(z)+\frac{1}{z}J'_m(z)+(1-\frac{m^2}{z^2})J_m(z)=0.
\end{equation}
Thus $J_m(z)$ has only simple zeroes.  On the other hand, by
multiplying both sides of \eqref{laurent} by $t^{-m-1}$ then
integrating on $|t|=1$, we obtain the following integral
representation of the Bessel functions:
\[
J_m(z)=\frac{1}{2\pi}\int_0^{2\pi}\cos(m\theta-z\sin\theta)\,
d\theta.
\]
From this integral representation, we know that $J_0(x)$ is
positive on the interval $[k\pi, (k+1/2)\pi]$ when $k$ is even and
negative on the interval when $k$ is odd.  It follows that
$J_0(x)$ has infinitely many of zeroes on the positive real axis
and all of these zeroes are on the intervals $((k+1/2)\pi,
(k+1)\pi)$. From \eqref{rec-2}, we know that
\begin{equation}\label{rec-3}
J_{m-1}(z)=z^{-m} \frac{d}{dz}(z^m J_m(z)),  \qquad
J_{m+1}(z)=-z^{m}\frac{d}{dz}(z^{-m}J_m(z)).
\end{equation}
It follows that $J_m(z)$ also has infinite many zeroes on the
positive real axis.  Furthermore, the zeroes of $J_m(z)$ and those
of $J_{m+1}(z)$ interlace.  Let $\lambda_{m, j}$, $j=1, 2,
\ldots$, be the positive zeroes of $J_m(z)$, arranged in
increasing order.  Then it follows from \eqref{bessel-eq} that
\[
\int_0^1 r J_m(\lambda_{m, j}r) J_m(\lambda_{m, k}r)\,
dr=\begin{cases} 0, \quad\quad\qquad\qquad j\not= k; \\
\frac{1}{2}J^2_{m+1}(\lambda_{m, j}), \quad j=k.
\end{cases}
\]
Furthermore, for any given integer $m$, $\{\sqrt{r}J_m(\lambda_{m,
j}r)\}_{j=1}^\infty$ forms a complete orthogonal basis for $L^2(0,
1)$.  Moreove, it follows from \eqref{rec-3} that for $m\ge 0$,
$\{r^{1/2+m}\}\cup\{\sqrt{r}J_m(\lambda_{m+1, j}r)\}_{j=1}^\infty$
forms a complete orthogonal basis for $L^2(0, 1)$ and so does
$\{\sqrt{r}J_m(\lambda_{m-1, j}r)\}_{j=1}^\infty$ for $m>0$.

\section{The computations}\label{comp}
Let $P=P(a_1, \ldots, a_n)=\{(z_1, \ldots, z_n)\in \C^n \mid
|z_1|<a_1, \ldots, |z_n|<a_n\}$. Write $\rho_j(z)=|z_j|^2-a_j^2$.
Then $P=\{z\in\C^n \mid \rho_j(z)<0, j=1, \ldots, n\}$.   Suppose
that
\[
u=\sumprime_{|J|=q} u_J d\bar z_J\in C^\infty(\ov{P}).
\]
For any integer $q$ between $1$ and $n-1$, we now solve the
$\dbar$-Neumann boundary value problem:
\begin{align}
\square_q &u=\lambda u; \label{lap}\\
&u\in \dom(\dbarstar_{q-1});\label{bdry-1}\\
\dbar_q &u\in \dom(\dbarstar_q).\label{bdry-2}
\end{align}
It follows from an easy integration by parts argument that $u\in
\dom(\dbarstar_{q-1})$ provided $u_{jK}(z)=0$ when $|z_j|=a_j$ for
any $(q-1)$-tuple $K$ and $j\in\{1, \ldots, n\}$. Write
$z_j=r_je^{i\theta_j}$. Using separation of variables, we write
\begin{equation}\label{sep}
u_J(z)=\prod_{k=1}^n u^k_J(z_k).
\end{equation}
Then $u\in\dom(\dbarstar_{q-1})$ provided
\begin{equation}\label{bdry-1a}
u^k_J(a_ke^{i\theta_k})=0, \qquad \text{when}\ k\in J.
\end{equation}
For any $K=(k_1, \ldots, k_{q+1})$, write
\[
v_K=\sum_{l=1}^{q+1} (-1)^{l+1} \frac{\partial u_{K\setminus
k_l}}{\partial \bar z_{k_l}},
\]
where $K\setminus k_l$ means the deletion of the $k_l$ entry from
$K$. Then
\[
\dbar u=\sumprime_{|K|=q+1} v_K d\bar z_K.
\]
Thus $\dbar_q u\in \dom(\dbarstar_q)$ if $v_{jJ}(z)=0$ whenever
$|z_j|=a_j$ for any $j\in\{1, \ldots, n\}$ and $q$-tuple $J$.
Using the separation of variables \eqref{sep}, we have that
$\dbar_q u\in \dom(\dbarstar_q)$ provided, in addition to
\eqref{bdry-1a}, $u_J$ also satisfies
\begin{equation}\label{bdry-2a}
\frac{\partial u^k_J}{\partial \bar z_k}(a_k e^{i\theta_k})=0,
\qquad \text{when}\ k\not\in J.
\end{equation}

Recall that $\square_q=(-1/4)\Delta$ where $\Delta$ is the usual
Laplacian acting diagonally. Denote by
$\Delta_k=4(\partial^2/\partial z_k\partial \bar z_k)$ the
Laplacian in the $z_k$-variable. Then, with the separation of
variables \eqref{sep}, the boundary value problem
\eqref{lap}-\eqref{bdry-2} is reduced to:
\begin{equation}\label{bdry-1b}
\Delta_k u^k_J(z_k)=-\lambda_k u^k_J, \qquad u^k_J(a_k
e^{i\theta_k})=0, \quad \text{for}\ k\in J,
\end{equation}
and
\begin{equation}\label{bdry-2b}
\Delta_k u^k_J(z_k)=-\lambda_k u^k_J(z_k), \qquad \frac{\partial
u^k_J}{\partial \bar z_k}(a_k e^{i\theta})=0, \quad \text{for}\
k\not\in J,
\end{equation}
with
\begin{equation}\label{lambda}
\lambda=\frac{1}{4}\sum_{k=1}^n \lambda_k.
\end{equation}
The boundary value problem \eqref{bdry-1b} gives the eigenvalues
for the Dirichlet Laplacian on the disc $|z_k|<a_k$. It is well
known (and easy to see) that these eigenvalues are
\begin{equation}\label{eigenv-1}
\left(\frac{\lambda_{m_k, j_k}}{a_k}\right)^2
\end{equation}
and the associated eigenfunctions are
\begin{equation}\label{eigenf-1}
J_{m_k}(\lambda_{m_k, j_k}r_k/a_k)e^{im_k\theta_k},
\end{equation}
for $m_k\in \Z$ and $j_k\in \N$.

To solve the boundary value problem \eqref{bdry-2b}, we separate
the variables in polar coordinates:
$u^k_J(z_k)=R(r_k)\Theta(\theta_k)$. Then \eqref{bdry-2b} is
reduced to
\begin{equation}\label{theta}
\frac{\Theta''}{\Theta}=-\mu, \qquad
\Theta(\theta_k+2\pi)=\Theta(\theta_k),
\end{equation}
and
\begin{equation}\label{r}
\frac{R''}{R}+\frac{1}{r_k}\frac{R'}{R}-\frac{\mu}{r^2_k}=-\lambda_k,\qquad
\frac{R'}{R}(a_k)=-\frac{i}{a_k}\frac{\Theta'}{\Theta}.
\end{equation}
From \eqref{theta}, we know that $\mu=m^2_k$, $m_k\in \Z$, with
the associated eigenfunctions $e^{im_k\theta_k}$.  We first
consider the case when $\lambda_k=0$. In this case, we know from
$\Theta=e^{im_k\theta_k}$ and \eqref{r} that $R=r_k^{m_k}$. Since
by interior elliptic regularity, the eigenfunctions must be smooth
at the origin, we know that $0$ is an eigenvalue of the boundary
value problem \eqref{bdry-2b} with the associated eigenfunctions
$z_k^{|m_k|}$.

Now we consider the case when $\lambda_k>0$. Using the
substitution $r=\sqrt{\lambda_k}r_k$, we reduce \eqref{r} to
\begin{equation}\label{r-1}
R''+\frac{1}{r}R'+(1-\frac{m_k^2}{r^2})R=0, \qquad
\sqrt{\lambda_k}a_kR'(\sqrt{\lambda_k}a_k)-m_kR(\sqrt{\lambda_k}a_k)=0.
\end{equation}
From \eqref{bessel-eq}, we know that $R=J_{m_k}(r)$, and from
\eqref{rec-2}, we know that $J_{m_k+1}(\sqrt{\lambda_k}a_k)=0$. In
summary, from the boundary value problem \eqref{bdry-2b}, we
obtain the eigenvalues
\begin{equation}\label{eigenv-2}
\left(\frac{\lambda_{m_k+1, j_k}}{a_k}\right)^2
\end{equation}
with the associated eigenfunctions
\begin{equation}\label{eigenf-2}
J_{m_k}(\lambda_{m_k+1, j_k}r_k/a_k)e^{im_k\theta_k},
\end{equation}
for $m_k\in \Z$ and $j_k\in \N$.

From the above computations, we now know that the spectrum of
$\square_q$ on the polydisc $P$ contains the eigenvalues
\begin{equation}\label{eigenv-3}
\frac{1}{4}\sum_{k\in J}\left(\frac{\lambda_{m_k,
j_k}}{a_k}\right)^2
\end{equation}
of infinite multiplicity with the associated eigenforms
\begin{equation}\label{eigenf-3}
\prod_{k\in J}\left(J_{m_k}(\lambda_{m_k,
j_k}r_k/a_k)e^{im_k\theta_k}\right) \prod_{k\not\in J}z_k^{|m_k|}
d\bar z_J,
\end{equation}
and eigenvalues
\begin{equation}\label{eigenv-4}
\frac{1}{4}\sum_{k=1}^n \left(\frac{\lambda_{m_k,
j_k}}{a_k}\right)^2
\end{equation}
with the associated eigenforms
\begin{equation}\label{eigenf-4}
\prod_{k\in J}\left(J_{m_k}(\lambda_{m_k, j_k}
r_k/a_k)e^{im_k\theta_k}\right)\prod_{k\not\in
J}\left(J_{m_k-1}(\lambda_{m_k,
j_k}r_k/a_k)e^{i(m_k-1)\theta_k}\right) d\bar z_J,
\end{equation}
for any strictly increase $q$-tuple $J$, $m_k\in \Z$, and
$j_k\in\N$.

It remains to show that the spectrum of $\square_q$ consists of
nothing else but the eigenvalues listed in \eqref{eigenv-3} and
\eqref{eigenv-4}. To do this, we use the following well known fact
from the general operator theory ({\it e.g.}, \cite{Davies95},
Lemma 1.2.2): Let $T$ be a symmetric operator on a complex Hilbert
space $H$. If there exists a complete orthonormal basis
$\{f_j\}_{j=1}^\infty$ and $\lambda_j\in \R$ such that
$Tf_j=\lambda_j f_j$, then $T$ is essentially self-adjoint and the
spectrum of $\ov{T}$ is the closure of
$\{\lambda_j\}_{j=1}^\infty$ in $\R$.  It follows from facts about
the Bessel functions stated in the last paragraph of
Section~\ref{prelim} that for each $q$-tuple $J$, the coefficients
of $d\bar z_J$ in \eqref{eigenf-3} and \eqref{eigenf-4} form a
complete orthogonal basis for $L^2(P)$. Thus the spectrum of
$\square_q$ contains nothing else but eigenvalues listed in
\eqref{eigenv-3} and \eqref{eigenv-4} with associated eigenforms
listed in \eqref{eigenf-3} and \eqref{eigenv-4} respectively.  The
bottom of the spectrum is
\[
\min_{|J|=q}\{\frac{\lambda^2_{0, 1}}{4}\sum_{k\in
J}\frac{1}{a_k^2}\},
\]
which is always of infinite multiplicity.

Since we now know explicitly the spectrum and the associated
eigenforms, it is then easy to explicitly express the
$\dbar$-Neumann operator as an infinite sum of projections onto
the eigenspaces.  We left this to the interested reader.

\bigskip

\noindent{\bf Acknowledgement:} The author thanks Professor Peter
Polyakov for stimulating discussions and kind encouragement.

\bibliography{survey}
% Bibliography generated by BibTeX with amsplain
% bibliographystyle and mrabbrev.bib abbreviations.
%
%%%

%% bibliography generated by BiBTeX
\providecommand{\bysame}{\leavevmode\hbox to3em{\hrulefill}\thinspace}

\end{document}